\documentclass{aip-cp-mod}

\usepackage[numbers]{natbib}
\usepackage{rotating}
\usepackage{graphicx}
\usepackage{cleveref}
\usepackage{color}

\definecolor{darkred}{RGB}{165,50,50}
\definecolor{darkgreen}{RGB}{0,100,0}

\newcommand\nD{\mathcal{D}}
\newcommand\nL{\mathcal{L}}
\newcommand\nS{\mathcal{S}}
\newcommand\nE{\mathcal{E}}
\def\ee{\mathrm{e}}

\begin{document}

% Title portion
\title{Efficient Adaptive Computation of the Dynamics of Mott Transistors}

\author[aff1]{Winfried Auzinger}
\eaddress{w.auzinger@tuwien.ac.at}
%\eaddress[url]{http://www.asc.tuwien.ac.at/~winfried}
\author[aff2]{Lukas Einramhof}
%\eaddress{e1525156@student.tuwien.ac.at}
\author[aff2]{Karsten Held}
\eaddress{held@ifp.tuwien.ac.at}
%\eaddress[url]{https://www.ifp.tuwien.ac.at/cms-1}
\author[aff2]{Anna Kauch}
\eaddress{kauch@ifp.tuwien.ac.at}
%\eaddress[url]{https://www.ifp.tuwien.ac.at/cms-1}
\author[aff3]{Othmar Koch}
\eaddress{othmar@othmar-koch.org}
%\eaddress[url]{http://www.othmar-koch.org}
\author[aff2]{Clemens Watzenb{\"o}ck}
\eaddress{clemens.watzenboeck@tuwien.ac.at}
%\eaddress[url]{https://www.ifp.tuwien.ac.at/cms-1}
\author[aff1]{Ewa Weinm{\"u}ller\corref{cor1}}
%\eaddress{e.weinmueller@tuwien.ac.at}
%\eaddress[url]{http://www.asc.tuwien.ac.at/~ewa}
\affil[aff1]{TU Wien, Institute of Analysis and Scientific Computing, Wiedner Hauptstrasse 8--10,
A-1040 Wien, Austria}
%Note the use of superscript ``a)'' to indicate the author's e-mail address below. Use b), c), etc. to %indicate e-mail addresses for more than 1 author.}
\affil[aff2]{TU Wien, Institute of Solid State Physics, Wiedner Hauptstrasse 8--10,
A-1040 Wien, Austria}
\affil[aff3]{University of Vienna, Wolfgang Pauli-Institute,
Oskar Morgenstern--Platz 1, A-1090 Wien, Austria}
\corresp[cor1]{Corresponding author: e.weinmueller@tuwien.ac.at}

\maketitle

\begin{abstract}
We investigate time-adaptive Magnus-type integrators for the numerical
approximation of a Mott transistor. The rapidly
attenuating electromagnetic field calls for adaptive choice of
the time steps. As a basis for step selection, asymptotically correct
defect-based estimators of the local error are employed. We analyze
the error of the numerical approximation in
the presence of the unsmooth external potential and demonstrate
the advantages of the adaptive approach.
\end{abstract}

% Head 1
\section{INTRODUCTION}
We study systems of linear ordinary differential equations of Schr{\"o}dinger type
\begin{equation}
\label{eq0}
\psi'(t) = -\,\mathrm{i}\,H(t)\,\psi(t)=:A(t)\psi(t)\,, \qquad t \in [t_0, t_{\mathrm{end}}]\,, \qquad
\psi(t_0) = \psi_0\; \mbox{ given}\,,
\end{equation}
with a time-dependent Hermitian matrix $H \colon \mathbb{R} \to \mathbb{C}^{d \times d}$.
The exact flow of (\ref{eq0}) is denoted by $\nE(t;t_0)\psi_0$ in the following.

We focus on a system that
is described by a Hubbard model of electron-hopping between discrete sites.
For the numerical experiments we use the parameters corresponding to a simple approximation of a Mott transistor~\cite{Zhong2015}.
The time dependence is introduced by coupling
of the electronic system to a source--drain potential which is switched on rapidly.
The corresponding electromagnetic field is treated classically and
enters the equation as a scalar potential and modifies hoppings via Peierls' substitution \cite{Li2020}.

We consider the numerical solution of the system by exponential-based
Magnus-type time integrators in conjunction with an adaptive Lanczos
method \cite{auzingeretal18b,jaweckietal20},
$$u_{n+1}=\nS(\tau;t_n)u_n, \qquad \mbox{where}$$
\begin{eqnarray*}
&&\nS(\tau;t_n) =
\nS_J(\tau) \cdots \nS_1(\tau)
= \ee^{\Omega_J (\tau)} \,\cdots\, \ee^{\Omega_1(\tau)}\,, \\
&&\Omega_j(\tau) = \tau\, B_j(\tau),~~ j=1, \ldots, J,
\qquad B_{j}(\tau) = \sum_{k=1}^K a_{jk}\, A_{k}(\tau),\qquad A_{k}(\tau) = A(t_n+c_k \tau)\,.
\end{eqnarray*}
The coefficients $a_{jk}$, $c_k$ are determined from the \emph{order conditions}
such that the method attains convergence order $ p $. In this study,
we will use the methods referred to as \texttt{CF2} (exponential midpoint
rule), \texttt{CF4oH} and \texttt{CF6n} in \cite{solar1}, the two
latter constructed by us with the aim to optimize the error for a given
computational effort. For comparison, we also show the embedded Runge--Kutta
method by Dormand \& Prince \texttt{DoPri45}.
The computational challenge results on the one hand from
the high dimension of the underlying system. Indeed, for a model with
$N$ discrete locations, the state space has dimension $4^N$. Thus, for
a model with the claim of physical relevance, the problem quickly reaches
the limitations of modern supercomputers.
On the other hand, the modeling of the switching process in a Mott transistor,
features a quickly attenuating electric field in a small time interval.
Thus the smoothness of the solution varies strongly in the course of the
time propagation, which suggests to use adaptive choice of the time step
to increase the efficiency and reliability of a numerical computation.
To this end, we use classical step-size choice based on asymptotically correct
defect-based estimators to equidistribute the local error, which were
constructed and analyzed in \cite{auzingeretal18b}.

\section{MODEL OF A MOTT TRANSISTOR}\label{sec:examples}

The Mott transistor, which is based on the Mott-insulator-to-metal transition driven by the change of the gate voltage $V_g$, can be described within the Hubbard model~\cite{hubbard63}, the paradigm model for the description of strongly
interacting electrons~\cite{Ma93,PKBS16}. We resort here to the \emph{second-quantization} formalism.
It describes the electron occupation on a given number of sites $N$, corresponding to Wannier discretization.
Only a single orbital per site is considered with nearest-neighbor hopping between the sites.
Due to the Pauli exclusion principle, there are only four states per site allowed (no electrons, one electron with spin-up or -down, two electrons with opposite spins). The electrons interact via Coulomb interaction $U_{ij\sigma\sigma'}$.
In the considered model the Hamiltonian in (\ref{eq0}) has the form
\begin{equation}\label{eq1}
H(t) = H_{\mathrm{stat}}+H_{\mathrm{dyn}}(t).
\end{equation}
The static part for arbitrary hopping  $v_{ij}$ reads
\begin{equation}\label{eq.HubHamStat}
	H_{\mathrm{stat}} = \sum_{i,j,\sigma} v_{ij} \hat{c}_{j\sigma}^{\dagger} \hat{c}_{i\sigma}^{\phantom\dagger}
	+ \frac{1}{2}\sum_{i,j,\sigma,\sigma'} U_{ij\sigma\sigma'} \hat n_{i\sigma} \hat n_{j\sigma'}
	%\hat{c}_{i\sigma}^{\dagger} \hat{c}_{j\sigma'}^{\dagger} \hat{c}_{j\sigma'} \hat{c}_{i\sigma}
    %,
\end{equation}
where $i,j$ sum over all $N$ sites  and the spins $\sigma,\sigma' \in\{\uparrow,\downarrow \}$
are either \emph{up} or \emph{down}.
The notation $\hat{c}_{j\sigma}^\dagger \hat{c}_{i\sigma}^{\phantom\dagger}$ describes the ``hopping'' of an electron
from site $i$ to $j$ with  \emph{creation} and \emph{annihilation operators}
$\hat{c}_{j \sigma}^\dagger$ and $\hat{c}_{i \sigma}$.
The \emph{hopping amplitudes} $v_{ij}$ with $i,j=1,\ldots,N$ give the probability (rate) of such an
electron hopping; $\hat{n}_{j\sigma}=\hat{c}_{j\sigma}^\dagger \hat{c}_{j\sigma}^{\phantom\dagger}$
is the \emph{occupation number operator}, which counts the number of electrons with spin $\sigma$ at site $j$.
For details on the notation in~(\ref{eq.HubHamStat}) refer to~\cite{hubbard63,Ma93,PKBS16}.
In the following we take only the local part of the Coulomb interaction, i.e. $U_{ij\sigma\sigma'}=U\delta_{ij}(1-\delta_{\sigma\sigma'})$.

\noindent The dynamic part of the Hamiltonian is in general given by
\begin{equation}\label{eq.HubHamDyn}
H_{\mathrm{dyn}}(t) = \sum_{\alpha=1}^{N_\alpha} \left(\Re (f_\alpha(t)) \:
H_{\mathrm{symm}}^{\alpha} + \mathrm{i} \Im (f_\alpha(t)) \: H_{\mathrm{anti}}^\alpha \right) +
\sum_{\beta=1}^{N_\beta} g_\beta(t) \: H_{\mathrm{pot}}^\beta.
\end{equation}
The real matrices $H_{\mathrm{symm}}^\alpha, H_{\mathrm{anti}}^\alpha,$ and
$H_{\mathrm{pot}}^\alpha$ have the following properties:
$ g_\alpha \in \mathbb{R}$, \
$ H_{\mathrm{symm}}^\alpha$ is symmetric, \
$ H_{\mathrm{anti}}^\alpha$ is skew-symmetric, \
$ H_{\mathrm{pot}}^\beta$ is diagonal.
Here, we restrict ourselves to the Coulomb gauge where
$f_\alpha(t) =0 \ \forall \alpha,$  $N_\beta=1.$
For the performance analysis we choose a rapidly attenuating potential
\begin{equation}
g_{1}(t) \equiv g(t) = 1 - \frac{1}{\mathrm{e}^{(t-t_0)/T}  +1},   \quad t_0,T \in \mathbb{R}^+,\qquad
H_{\mathrm{pot}}^1 = \sum_{i\sigma} V_i \hat{n}_{i \sigma} =: H_{\mathrm{pot}},
\quad V_i \in \mathbb{R}.\label{pulse1}
\end{equation}

\section{ERROR ANALYSIS}\label{sec:ana}

In this section, we investigate the error structure of commutator-free
Magnus-type time integrators for the model of a Mott transistor.
To this end, we recapitulate general convergence results given in \cite{auzingeretal18b}:

The local time-stepping error
can be expressed in terms of an iterated defect $\nD_i(\tau;t_0)$, where
$$ \nD_0=\nD=\psi'(t)-A(t)\psi(t), \qquad \nD_{i+1} = \nD_i'-A\nD_i,\quad i=0,1,\dots.$$
For the exponential midpoint rule, we thereby obtain for the local error
$\nL(\tau;t_0)\psi_0:=(\nS(\tau;t_0)-\nE(\tau;t_0))\psi_0$
\begin{eqnarray*}
 \nL(\tau;t_0)\hspace*{-3mm}&=& \hspace*{-3mm}\int_{0}^{\tau} \Pi(\tau,\sigma_1)\,\nD(\sigma_1;t_0)\,\mathrm{d}\sigma_1 =
\int_{0}^{\tau} \Pi(\tau,\sigma_1)\,\int_0^{\sigma_1}\Pi(\sigma_1,\sigma_2)\,
\int_0^{\sigma_2}\Pi(\sigma_2,\sigma_3)\, \mathrm{d}\sigma_3\, \mathrm{d}\sigma_2\,\mathrm{d}\sigma_1
 \cdot \nD_2(0;t_0) + O(\tau^4)\\
\hspace*{-3mm}&=:& \hspace*{-1mm}\underbrace{\mathcal{I}_3(\tau)}_{=\,O(\tau^3)}\!\cdot\,(\Gamma''(0)-A''(t_0))
+ O(\tau^4),\qquad \mbox{where} \quad\frac{\mathrm{d}}{\mathrm{d}\tau} \mathrm{e}^{\tau\, A(\tau/2)}
=\Gamma(\tau)\,\mathrm{e}^{\tau\,A(\tau/2)}\,.
\end{eqnarray*}
It follows from the arguments in \cite{auzingeretal18b} that more generally,
the leading local error term for a higher-order CFM of order $p$ depends
on the difference $\Gamma^{(p)}(0)-A^{(p)}(0)$.

For the exponential midpoint rule of order two, this means in particular that the following
error estimate holds:\\[2mm]
\cite[Proposition 4.1]{auzingeretal18b} \label{pro:EMR-locerr}
\textit{Consider the solution of~{\rm (\ref{eq0})} by the exponential
midpoint rule.
If $A\in C^3$, then the local error satisfies
\begin{equation} \label{EMR-L-leading}
\| \nL(\tau;t_0) \|_2 \leq
\frac{1}{12} \tau^3 \big\| [A(t_0)),A'(t_0)]
                         - \frac{1}{2} A''(t_0) \big\|_2
+ O(\tau^4)\,.
\end{equation}
}
For the Hubbard model with the potential (\ref{pulse1}), the convergence analysis above
has the following implications.

Differentiation of $g(t)$ shows the following asymptotics:
\begin{eqnarray}
g(t)       &=& 1- \frac{1}{\mathrm{e}^{(t-t_0)/T}+1} \in (0,1)\,,\\
g'(t)      &=& \frac{1}{T} \frac{\mathrm{e}^{(t-t_0)/T}}{(\mathrm{e}^{(t-t_0)/T}+1)^2} \in (0,1/T),\qquad
g''(t)     = \frac{1}{T^2} \left(\frac{\mathrm{e}^{(t-t_0)/T}}{(\mathrm{e}^{(t-t_0)/T}+1)^2} -
                 \frac{2 \mathrm{e}^{(t-t_0)/T}}{(\mathrm{e}^{(t-t_0)/T}+1)^3}\right) \in (0,1/T^2)\,,\\
%g^{(3)}(t) &=& \frac{1}{T^3} \left(\frac{\mathrm{e}^{(t-t_0)/T}}{(\mathrm{e}^{(t-t_0)/T}+1)^2} - \frac{6\mathrm{e}^{2(t-t_0)/T}}{(\mathrm{e}^{(t-t_0)/T}+1)^3}
%                 +6\frac{\mathrm{e}^{3(t-t_0)/T}}{(\mathrm{e}^{(t-t_0)/T}+1)^4}\right) \\
             &\vdots& \nonumber \\
g^{(p)}(t) &=&  O(T^{-p})\,.
\end{eqnarray}
When considering the asymptotics of the error as a function of $T$, clearly the local
error is dominated by the highest derivative of $A(t)$, where the time-dependence
occurs in $g$. For the exponential midpoint
rule, for example, (\ref{EMR-L-leading}) shows that this is $A''$. Analogously,
for a method of order $p$ the local error
is dominated by $\Gamma^{(p)}(0)-A^{(p)}(0)$ and is thus proportional to $\tau^{p+1}T^{-p}$.

\section{NUMERICAL RESULTS}\label{sec:num}

\begin{center}
\begin{figure}
\includegraphics[height=1.4cm, trim=0mm -45mm 0mm 45mm]{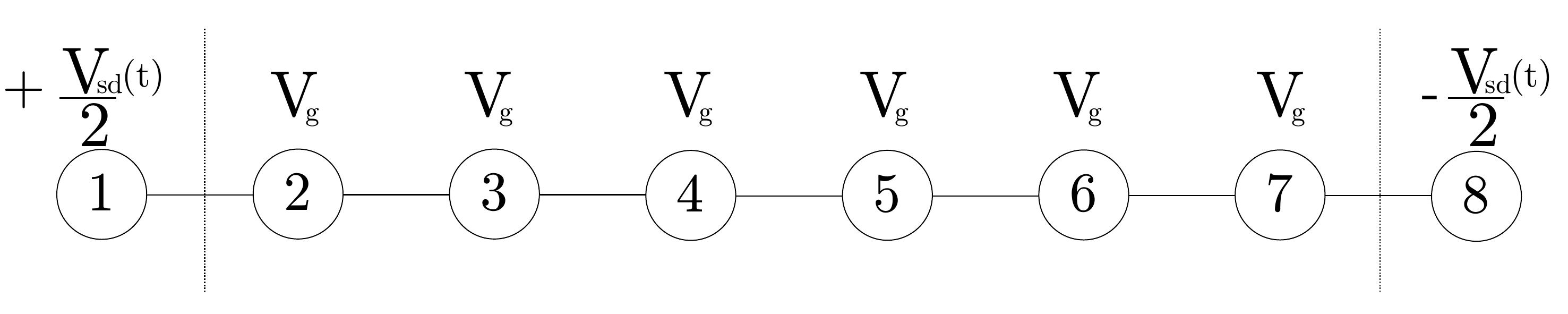} \hspace*{2mm}
\includegraphics[height=3.0cm, trim=0mm 3mm 0mm 0mm,clip]{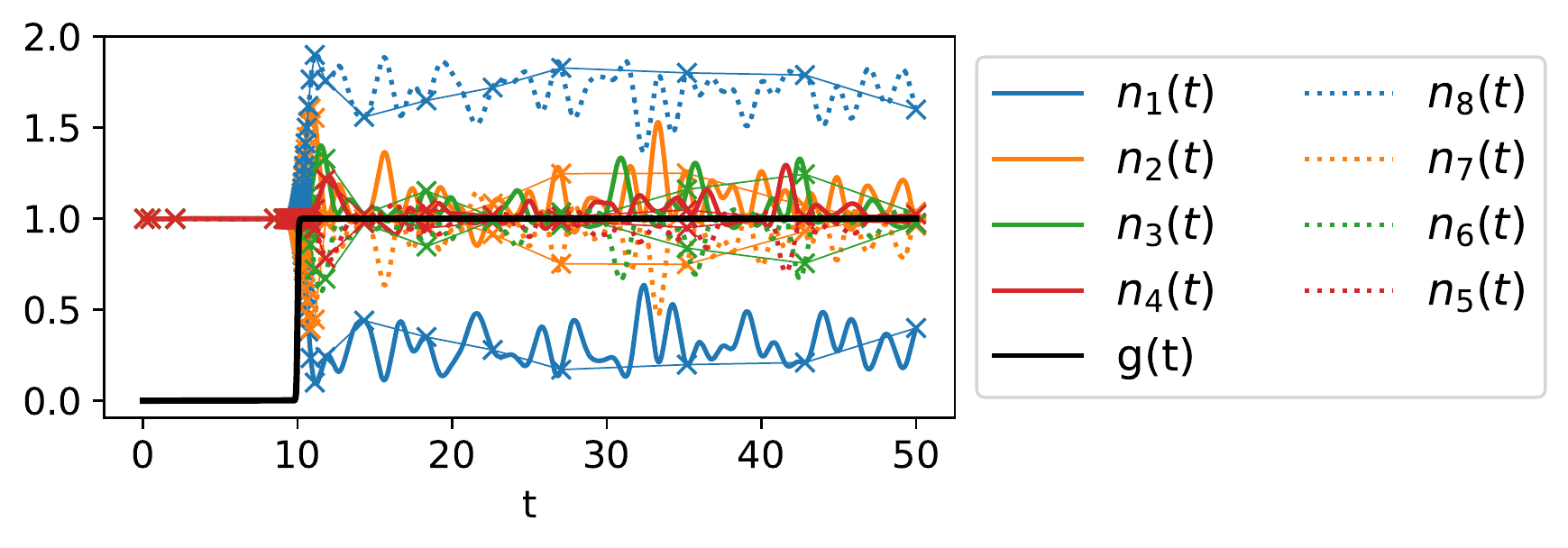}
\caption{\emph{Left}: Schematic depiction of the Mott transistor with the source--drain voltage applied to the outer sites [described by the potential $V_{sd}(t)$]  and gate voltage ($V_g$).
\emph{Right}: Time dependence of the charge  at the sites $n_i(t)$ and of the source--drain potential $g(t)$.
Computed with $\texttt{DoPri45}$ (solid lines) and $\texttt{CF4oH}$ (crosses).}
\label{fig:model}
\end{figure}
\end{center}
\begin{center}
\begin{figure}
\includegraphics[height=3.5cm]{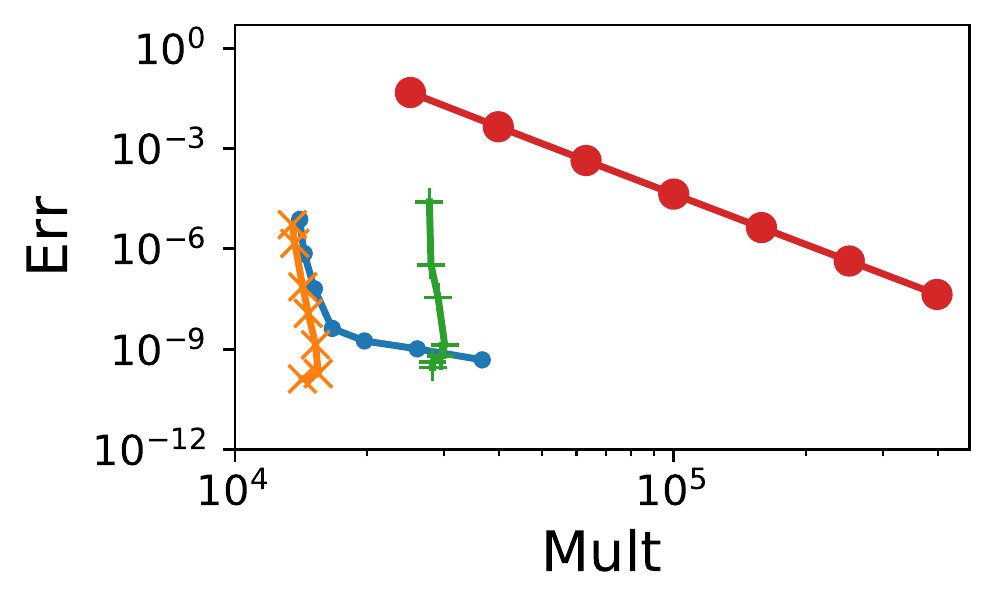} \hspace*{0.5cm}
\includegraphics[height=3.5cm]{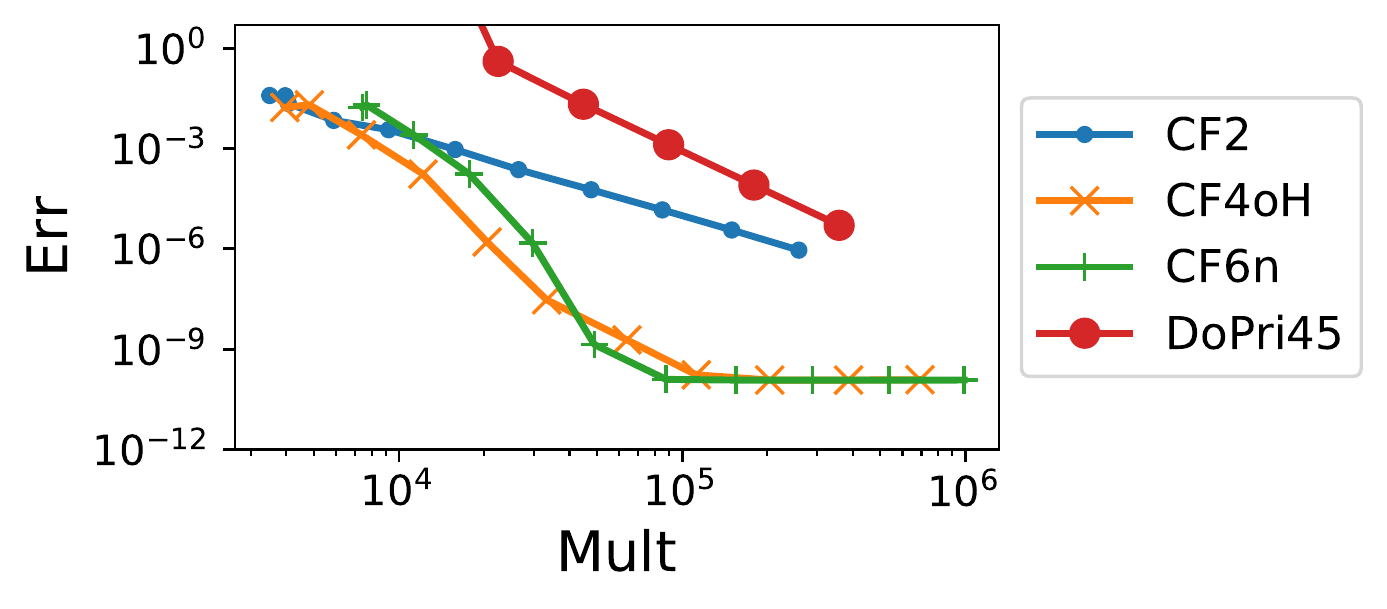}
\caption{Comparison of the error as a function of the number of matrix--vector multiplications
for adaptive time stepping (left) and equidistant time steps (right).}
\label{fig:comp}
\end{figure}
\end{center}

As a numerical illustration, we show the adaptive solution of our model of a Mott transistor
by Magnus-type integrators. The geometry of the model is illustrated in Figure~\ref{fig:model} (left).
The dynamics of the transistor is approximated by the dynamics of an 8-site  chain with 8 electrons.
%This defines the static part
%$H_\mathrm{stat}$ of the Hamiltonian (\ref{eq1}).
The sites are then split into three sections which correspond to an electron reservoir (left and right) which is coupled to a source--drain potential $V_{sd}$. The scalar potential added on the mid-sites ($V_g$ in the left part of Figure~\ref{fig:model}) models the gate of the transistor. It controls whether  the transistor is in a conducting or an insulating state. For the numerical analysis shown here we chose $V_g=0=V_i, \quad i \in {2,3,4,5,6,7}$ (conducting state). For the source--drain potential we used $V_1=-V_8=V_{sd}/2=1.04 U$; $T=2^{-5}$. The potential is switched on at $t_0=10$ and $\psi_0$ is chosen as the (unique) ground state of the model. Hubbard interaction $U=10 |v_{12}|$.
The error displayed is
$||{\psi}(t=50) - {\psi}^{\mathrm{ref}}(t=50) ||_2$,
where the reference solution $\psi^{\mathrm{ref}}$ was computed
with adaptive \texttt{CF4oH} with a high accuracy of $\texttt{tol}=10^{-11}$.
The tolerance of the Lanczos iteration was $10^{-12}$ throughout.
In the right plot of Figure~\ref{fig:model}, the index of $n_i=n_{i\uparrow}+n_{i\downarrow}$
corresponds to the sites labeled in the left part of the figure, and
$n_{i\sigma}(t) = \psi^\dagger(t) \, \hat n_{i\sigma} \, \psi(t)$.
It shows the occupation numbers (proportional to the local charge) computed
with a fine error tolerance by the Dormand \& Prince Runge--Kutta method
(solid lines), and the adaptive Magnus-type method \texttt{CF4oH} (crosses)
with the same error tolerance of $\mathtt{tol}=10^{-11}$.
We observe that the solution computed with adaptive time stepping
chooses appropriately small steps where the solution varies strongly,
and very large steps where the solution is smooth. Nonetheless,
the solutions correspond very well at the grid points. We conclude
that adaptive commutator-free Magnus-type methods serve very well
to accurately approximate the solution both in regions where dense
grids are required and where the dynamics allows for large time steps.
The adaptive choice of the steps thus implies a gain in efficiency
in contrast to uniform time steps.
Indeed, in Figure~\ref{fig:comp}, we show the error as a function of the
number of matrix--vector multiplications, which constitute most of the
computational effort required for the integrators.
The left plot gives the results for adaptive time stepping, on
the right-hand side, the results for equidistant time steps are
given for comparison. We observe that the error
is smaller for a comparable computational effort when the adaptive
strategy is used. This results from the unsmooth solution dynamics
associated with the rapid changes in the electric field. The specially constructed
CFM integrators are most efficient, adaptive Runge--Kutta methods
show consistent convergence behavior, but prohibitively large computational effort.

%\begin{center}
%\begin{figure}
%\includegraphics[width=5cm]{8x1_occ_iT5.pdf}
%\caption{\ask{What is $n_i$? What is the meaning of the index? Anyway, the plot looks great!}.\label{fig:steps}}
%\end{figure}
%\end{center}

\section*{ACKNOWLEDGMENTS}
This work was supported by the Austrian Science Fund (FWF) under grant P 30819-N32.
%The computations have been performed on the Vienna Scientific Cluster (VSC).

%\bibliographystyle{aipnum-cp}%
%\bibliography{num,schroedinger,books,appl,gen,held,KK,karsten-key2,watzenboeck}%

%Control: year  (1) truncated

\end{document}